\RequirePackage{fix-cm}
\documentclass[smallextended, reqno, 12pt]{svjour3}       
\usepackage{amssymb,amsmath,amscd, amsfonts,latexsym}
\usepackage[left=4cm, right=4cm]{geometry}
\usepackage{graphicx}
\usepackage{url}
\usepackage{hyperref}
\usepackage{float}

\begin{document}

\title{Repdigits as sums of three Padovan numbers}

\thanks{The author was supported by the Austrian Science Fund (FWF) grants: F5510-N26 -- Part of the special research program (SFB), ``Quasi Monte Carlo Metods: Theory and Applications'' and W1230 --``Doctoral Program Discrete Mathematics''.}


\titlerunning{Repdigits as sums of three Padovan numbers}        

\author{Mahadi Ddamulira}


\institute{M. Ddamulira \at
              Institute of Analysis and Number Theory, Graz University of Technology, Kopernikusgasse 24/II, A-8010 Graz, Austria. \\
              \email{mddamulira@tugraz.at; mahadi@aims.edu.gh}           
}

\date{Received: date / Accepted: date}

\maketitle

\begin{abstract}
Let $ \{P_{n}\}_{n\geq 0} $ be the sequence of Padovan numbers  defined by $ P_0=0 $, $ P_1 =1=P_2$ and $ P_{n+3}= P_{n+1} +P_n$ for all $ n\geq 0 $. In this paper, we find all repdigits in base $ 10 $ which can be written as a sum of three Padovan numbers. 
\end{abstract}

\section{Introduction}
Let $ \{P_{n}\}_{n\geq 0} $ be the sequence of  Padovan numbers given by
$$P_0=0, ~P_1= 1, ~P_2=1 \text{  and  } P_{n+3}= P_{n+1}+P_n \text{   for all  } n\geq 0.$$
The first few terms of this sequence are
$$\{P_{n}\}_{n\ge 0} = 0, 1, 1, 1, 2, 2, 3, 4, 5, 7, 9, 12, 16, 21, 28, 37, 49, 65, 86, 114, 151, \ldots. $$

A repdigit is a positive integer $ N $ that has only one distinct digit when written in base $ 10 $. That is, $ N $ is of the form
\begin{eqnarray}\label{rep1}
N=d\left(\dfrac{10^{\ell}-1}{9}\right)
\end{eqnarray}
for some positive integers $ d, \ell $ with $ 1\le d\le 9 $.

\section{Main Result}
In this paper, we study the problem of writing repdigits as sums of three Balancing numbers. More prcisely, we completely solve the Diophantine equation
\begin{eqnarray}
N=P_{n_1}+P_{n_2}+P_{n_3}=d\left(\dfrac{10^{\ell}-1}{9}\right),\label{Problem}
\end{eqnarray}
in non-negative integers $ (N, n_1, n_2, n_3, d, \ell ) $ with $ n_1\ge n_2\ge n_3 \ge 0 $, $ \ell\ge 2 $ and $ 1\le d\le 9 $.

We discard the situations when $ n_1=1 $ and $ n_1=2 $ and just count the solutions for $ n_1=3 $ since $ P_1=P_2=P_3 =1 $. For the same reasons, we discard the situation when $ n_1=4 $ and just count the solutions for $ n_1=5 $ since $ P_4=P_5 =2 $.  Thus, we always assume that $ n_1,n_2, n_3 \notin \{1,2,4\}$. Our main result is the following.
\begin{theorem}\label{Main}
All non-negative interger solutions $ (N, n_1, n_2, n_3, d, \ell ) $ with $ n_1\ge n_2\ge n_3 \ge 0 $, $ \ell\ge 2 $ and $ 1\le d\le 9 $  to the Diophantine equation \eqref{Problem} arise from
$$ N\in\{11, 22, 33, 44, 55, 66, 77, 88, 99, 111, 222, 333, 444, 555, 666, \
888, 1111, 3333, 7777\}.$$
\end{theorem}

This paper serves as a continuation of the results in \cite{Ddamulira1}, \cite{Lomeli}, \cite{Luca1}, \cite{Luca2}, \cite{Luca3} and \cite{Luca4}. The method of proof involves the application of Baker's theory for linear forms in logarithms of algebraic numbers, and the Baker-Davenport reduction procedure. Computations are done using a simple computer program in \textit{Mathematica}.

\section{Preliminary results}
\subsection{The Padovan sequence} 
Here, we recall some important properties of the Padovan sequence $ \{P_n\}_{n\geq 0} $. The characteristic equation
\begin{eqnarray*}
 x^3-x-1 = 0
\end{eqnarray*}
has roots $ \alpha, \beta, \gamma = \bar{\beta} $, where
\begin{eqnarray}\label{Pado1}
\alpha =\dfrac{r_1+r_2}{6}, \qquad \beta = \dfrac{-(r_1+r_2)+\sqrt{-3}(r_1-r_2)}{12}
\end{eqnarray}
and
\begin{eqnarray}\label{Pado2}
r_1=\sqrt[3]{108+12\sqrt{69}} \quad \text{and}\quad r_2=\sqrt[3]{108-12\sqrt{69}}.
\end{eqnarray}
Furthermore, the Binet formula is given by
\begin{eqnarray}\label{Pado3}
P_n = a\alpha^{n}+b\beta^{n}+c\gamma^{n} \qquad \text{ for all} \quad n\ge 0,
\end{eqnarray}
where
\begin{eqnarray}\label{Pado4}
\quad a=\dfrac{\alpha+1}{(\alpha-\beta)(\alpha-\gamma)}, \quad b= \dfrac{\beta +1}{(\beta -\alpha)(\beta-\gamma)}, \quad c = \dfrac{\gamma+1}{(\gamma-\alpha)(\gamma-\beta)}=\bar{b}.
\end{eqnarray}
Numerically, the following estimates hold:
\begin{eqnarray}\label{Pado5}
&1.32<\alpha<1.33\nonumber\\
&0.86 < |\beta|=|\gamma|=\alpha^{-\frac{1}{2}}< 0.87\\
&0.72<a<0.73\nonumber\\
&0.24<|b|=|c|<0.25.\nonumber
\end{eqnarray}
From \eqref{Pado1}, \eqref{Pado2} and \eqref{Pado5}, it is easy to see that the contribution the complex conjugate roots $ \beta $ and $ \gamma $, to the right-hand side of \eqref{Pado3}, is very small. In particular, setting
\begin{eqnarray}\label{Pado6}
e(n):=P_n-a\alpha^{n}=b\beta^{n}+c\gamma^{n}\quad \text{ then } \quad |e(n)|< \dfrac{1}{\alpha^{n/2}}
\end{eqnarray}
holds for all $ n\ge 1 $.
Furthermore, by induction, we can prove that 
\begin{eqnarray}\label{Pado7}
\alpha^{n-3}\leq P_n \leq \alpha^{n-1} \quad \text{holds for all }\quad n\geq 1.
\end{eqnarray}

\subsection{Linear forms in logarithms}
Let $ \eta $ be an algebraic number of degree $ D $ with minimal primitive polynomial over the integers
$$ a_{0}x^{d}+ a_{1}x^{d-1}+\cdots+a_{d} = a_{0}\prod_{i=1}^{D}(x-\eta^{(i)}),$$
where the leading coefficient $ a_{0} $ is positive and the $ \eta^{(i)} $'s are the conjugates of $ \eta $. Then the \textit{logarithmic height} of $ \eta $ is given by
$$ h(\eta) := \dfrac{1}{D}\left( \log a_{0} + \sum_{i=1}^{d}\log\left(\max\{|\eta^{(i)}|, 1\}\right)\right).$$

In particular, if $ \eta = p/q $ is a rational number with $ \gcd (p,q) = 1 $ and $ q>0 $, then $ h(\eta) = \log\max\{|p|, q\} $. The following are some of the properties of the logarithmic height function $ h(\cdot) $, which will be used in the next sections of this paper without reference:
\begin{eqnarray}
h(\eta_1\pm \eta_2) &\leq& h(\eta_1) +h(\eta_2) +\log 2,\nonumber\\
h(\eta_1\eta_2^{\pm 1})&\leq & h(\eta_1) + h(\eta_2),\\
h(\eta^{s}) &=& |s|h(\eta) ~~~~~~ (s\in\mathbb{Z}). \nonumber
\end{eqnarray}

\begin{theorem}\label{Matveev11} Let $\eta_1,\ldots,\eta_t$ be positive real algebraic numbers in a real algebraic number field  $\mathbb{K} \subset \mathbb{R}$ of degree $D_\mathbb{K}$, $b_1,\ldots,b_t$ be nonzero integers, and assume that
\begin{equation}
\label{eq:Lambda}
\Lambda:=\eta_1^{b_1}\cdots\eta_t^{b_t} - 1,
\end{equation}
is nonzero. Then
$$
\log |\Lambda| > -1.4\times 30^{t+3}\times t^{4.5}\times D_\mathbb{K}^{2}(1+\log D_{\mathbb{K}})(1+\log B)A_1\cdots A_t,
$$
where
$$
B\geq\max\{|b_1|, \ldots, |b_t|\},
$$
and
$$A
_i \geq \max\{D_{\mathbb{K}} h(\eta_i), |\log\eta_i|, 0.16\},\qquad {\text{for all}}\qquad i=1,\ldots,t.
$$
\end{theorem}

\subsection{Reduction procedure}\label{Reduction}
During the calculations, we get upper bounds on our variables which are too large, thus we need to reduce them. To do so, we use some results from the theory of continued fractions. 

For the treatment of linear forms homogeneous in two integer variables, we use the well-known classical result in the theory of Diophantine approximation.
\begin{lemma}\label{Legendre}
Let $\tau$ be an irrational number,  $ \frac{p_0}{q_0}, \frac{p_1}{q_1}, \frac{p_2}{q_2}, \ldots $ be all the convergents of the continued fraction of $ \tau $ and $ M $ be a positive integer. Let $ N $ be a nonnegative integer such that $ q_N> M $. Then putting $ a(M):=\max\{a_{i}: i=0, 1, 2, \ldots, N\} $, the inequality
\begin{eqnarray*}
\left|\tau - \dfrac{r}{s}\right|> \dfrac{1}{(a(M)+2)s^{2}},
\end{eqnarray*}
holds for all pairs $ (r,s) $ of positive integers with $ 0<s<M $.
\end{lemma}

For a nonhomogeneous linear form in two integer variables, we use a slight variation of a result due to Dujella and Peth{\H o} (see \cite{dujella98}, Lemma 5a). For a real number $X$, we write  $||X||:= \min\{|X-n|: n\in\mathbb{Z}\}$ for the distance from $X$ to the nearest integer.
\begin{lemma}\label{Dujjella}
Let $M$ be a positive integer, $\frac{p}{q}$ be a convergent of the continued fraction of the irrational number $\tau$ such that $q>6M$, and  $A,B,\mu$ be some real numbers with $A>0$ and $B>1$. Let further 
$\varepsilon: = ||\mu q||-M||\tau q||$. If $ \varepsilon > 0 $, then there is no solution to the inequality
$$
0<|u\tau-v+\mu|<AB^{-w},
$$
in positive integers $u,v$ and $w$ with
$$ 
u\le M \quad {\text{and}}\quad w\ge \dfrac{\log(Aq/\varepsilon)}{\log B}.
$$
\end{lemma}

Finally, the following Lemma is also useful. It is Lemma 7 in \cite{guzmanluca}. 
\begin{lemma}
\label{gl}
If $r\geqslant 1$, $H>(4r^2)^r$  and $H>L/(\log L)^r$, then
$$
L<2^rH(\log H)^r.
$$
\end{lemma}

\section{Bounding the variables}
We assume that $ n_1\ge n_2\ge n_3 $. From \eqref{Problem} and \eqref{Pado7}, we have
\begin{eqnarray*}
\alpha^{n_1-3}\le P_{n_1} \le P_{n_1}+P_{n_2}+P_{n_3} = d\left(\dfrac{10^{\ell}-1}{9}\right) \le 10^{\ell}
\end{eqnarray*}
and 
\begin{eqnarray*}
10^{\ell-1} \le d\left(\dfrac{10^{\ell}-1}{9}\right)  =  P_{n_1}+P_{n_2}+P_{n_3} \le 3P_{n_1} < \alpha^{n+3},
\end{eqnarray*}
where we use $ \alpha^4 >3 $. Thus,
\begin{eqnarray*}
(n_1-3)\dfrac{\log\gamma}{\log 10} \le \ell \quad \text{and}\quad \ell -1 \le (n_1+3)\dfrac{\log\gamma}{\log 10}.
\end{eqnarray*}
Since $ \log \gamma /\log 10 = 0.122123... <1/5 $, we can conclude from the above that
\begin{eqnarray}\label{bal3}
(n_1-3)/5 < \ell < (n_1+8)/5.
\end{eqnarray}
Running a \textit{Mathematica} program in the range $ 0\le n_3 \le n_2 \le n_3 \le 500 $, $ 1 \le d \le 9 $ and $ 1 \le \ell \le 100 $ we obtain only the solutions listed in Theorem \ref{Main}. From now onwards, we assume that $ n_1> 500 $.

By using \eqref{Pado6}, equation \eqref{Problem} can be written as
\begin{eqnarray}\label{bal4}
a\alpha^{n_1}+e(n_1)+a\alpha^{n_2}+e(n_2)+a\alpha^{n_3}+e(n_3) =d \left(\dfrac{ 10^{\ell}-1}{9}\right).
\end{eqnarray}
We then consider \eqref{bal4} in three different cases as follows.
\subsection{Case 1} We have that
\begin{eqnarray*}
a\alpha^{n_1}+e(n_1)+a\alpha^{n_2}+e(n_2)+a\alpha^{n_3}+e(n_3) -\dfrac{d\cdot 10^{\ell}}{9} = -\dfrac{d}{9}.
\end{eqnarray*}
This is equivalent to 
\begin{eqnarray*}
a\alpha^{n_1}- \dfrac{d\cdot 10^{\ell}}{9} = -\dfrac{d}{9} -a(\alpha^{n_2}+\alpha^{n_3})- e(n_1)-e(n_2)-e(n_3).
\end{eqnarray*}
Thus, we have 
\begin{eqnarray*}
\left|a\alpha^{n_1}- \dfrac{d\cdot 10^{\ell}}{9}\right| &\le  &\dfrac{d}{9}+ a(\alpha^{n_2}+\alpha^{n_3})+|e(n_1)|+|e(n_2)|+|e(n_3)|\\
&<& 1+2a\alpha^{n_2}+3\alpha^{-n_3/2}\\
&<& 5a\alpha^{n_2},
\end{eqnarray*}
and so
\begin{eqnarray}\label{bal5}
\left|a\alpha^{n_1}- \dfrac{d\cdot 10^{\ell}}{9}\right|  < 5a\alpha^{n_2}.
\end{eqnarray}
We divide through \eqref{bal5} by $ a\alpha^{n_1}$ to get
\begin{eqnarray*}
\left| 10^{\ell}\alpha^{-n_1}\left(\dfrac{d}{9a}\right)-1\right| &<& 5\alpha^{n_2-n_1}. 
\end{eqnarray*}
Thus, we have
\begin{eqnarray}\label{bal6}
\left| 10^{\ell}\alpha^{-n_1}\left(\dfrac{d}{9a}\right)-1\right|< \dfrac{5}{\alpha^{n_1-n_2}}.
\end{eqnarray}
We put
\begin{eqnarray*}
\Lambda_1:= 10^{\ell}\alpha^{-n_1}\left(\dfrac{d}{9a}\right)-1.
\end{eqnarray*}
In order to apply Theorem \ref{Matveev11} we need to check that $ \Lambda_1\ne 0 $. Suppose that $ \Lambda_1=0 $, then we have
\begin{eqnarray}\label{gal1}
a\alpha^{n_1} = \dfrac{10^{\ell}\cdot d}{9}.
\end{eqnarray}
To see that this is not true, we consider the $ \mathbb{Q} $-automorphism $ \sigma $ of the Galois extension $ \mathbb{Q}(\alpha, \beta) $ over $ \mathbb{Q} $ given by $ \sigma(\alpha):=\beta $ and $ \sigma(\beta):=\alpha $. Now, if $ \Lambda_1=0 $, then $ \sigma(\Lambda_1)=0 $.
Thus, conjugating the relation \eqref{gal1}  under $ \sigma $, and taking absolute values on both sides, we get 
\begin{eqnarray*}
\dfrac{10^{\ell}\cdot d}{9} = |\sigma (a\alpha^{n_1})|=|b||\beta|^{n_1}<  |b|< \dfrac{1}{3},
\end{eqnarray*}
which is false for $ \ell \ge 2 $ and $ d \ge 1 $. Therefore, $ \Lambda_1\ne 0 $.

So we apply Theorem \ref{Matveev11} with the data
\begin{eqnarray*}
&t:=3, \quad \eta_1:=10, \quad  \eta_2:=\alpha, \quad \eta_3:=\dfrac{d}{9a},\quad  b_1:=\ell, \quad b_2:=-n_1, \quad b_3:=1.
\end{eqnarray*}
It is a well--known fact that  $$ a=\dfrac{\alpha(\alpha+1)}{3\alpha^2-1}, $$ the mimimal polynomial of $ a $ is $ 23x^3-23x^2+6x-1 $ and has roots $ a, b, c $. Since $  |b|=|c|< |a|=a <1 $ (by \eqref{Pado5}), then
\begin{eqnarray*}
h(a)=\dfrac{1}{3}\log 23.
\end{eqnarray*}
Since $ \eta_1, \eta_2, \eta_2\in  \mathbb{Q}(\alpha) $, we take the field $ \mathbb{K}:= \mathbb{Q}(\alpha) $ with degree $ D_{\mathbb{K}}:=3 $. Since $ \max\{1, \ell, n_1\} \le n_1$, we take $ B:=n_1 $. Further, the minimal polynomial of $ \alpha $ over $ \mathbb{Z} $ is $ x^3-x-1 $ has roots $ \alpha, ~ \beta, ~ \gamma $  with $ 1.32 < \alpha < 1.33 $ and $ |\beta|=|\gamma| <1 $. Thus, we can take $  h(\alpha) = \frac{1}{3} \log \alpha $. Similarly, $ h(10) = \log 10 $ . Also,
\begin{eqnarray*}
h(\eta_3)&\le& h(d)+h(9) + h(a) \le 4\log 3+\dfrac{1}{3}\log 23 < 5\log 3
\end{eqnarray*}
Thus, we can take $ A_1:= 3\log 10 $, $ A_2:=\log \alpha $ and $ A_3:=15\log 3 $. So, Theorem \ref{Matveev11} tells us that the left-hand side of \eqref{bal6} is bounded below by
\begin{eqnarray*}
\log |\Lambda_1|&>&-1.4 \times 30^{6}\times 3^{4.5}\times 3^{2}(1+\log 3)(1+\log n_1)(3\log 10)(\log \alpha)(15\log 3)\\
&>&-6.16\times 10^{14}\log n_1 \log\alpha.
\end{eqnarray*}
By comparing the above inequality with the right-hand side of \eqref{bal6} we get that 
\begin{eqnarray}\label{bala1}
n_1-n_2 \le 6.18 \times 10^{14}\log n_1.
\end{eqnarray}
\subsection{Case 2}
 We have that
\begin{eqnarray*}
a\alpha^{n_1}+e(n_1)+a\alpha^{n_2}+e(n_2) -\dfrac{d\cdot 10^{\ell}}{9} = -\dfrac{d}{9}-a\alpha^{n_3}-e(n_3).
\end{eqnarray*}
This is equivalent to 
\begin{eqnarray*}
a(\alpha^{n_1} +\alpha^{n_2})- \dfrac{d\cdot 10^{\ell}}{9} = -\dfrac{d}{9} -a\alpha^{n_3}- e(n_1)-e(n_2)-e(n_3).
\end{eqnarray*}
Thus, we have 
\begin{eqnarray*}
\left|a(\alpha^{n_1}+\alpha^{n_2})- \dfrac{d\cdot 10^{\ell}}{9}\right| &\le  &\dfrac{d}{9}+ a\alpha^{n_3}+|e(n_1)|+|e(n_2)|+|e(n_3)|\\
&<& 1+a\alpha^{n_3}+3\alpha^{-n_3/2}\\
&<& 3a\alpha^{n_3},
\end{eqnarray*}
and so
\begin{eqnarray}\label{bal8}
\left|a(\alpha^{n_1}+\alpha^{n_2})- \dfrac{d\cdot 10^{\ell}}{9}\right|  < 3a\alpha^{n_3}.
\end{eqnarray}
We divide through \eqref{bal5} by $ a(\alpha^{n_1}+\alpha^{n_2})$ to get
\begin{eqnarray*}
\left| 10^{\ell}\alpha^{-n_2}\left(\dfrac{d}{9a(1+\alpha^{n_1-n_2})}\right)-1\right| &<& \dfrac{3\alpha^{n_3-n_2}}{1+\alpha^{n_1-n_2}}. 
\end{eqnarray*}
Thus, we have
\begin{eqnarray}\label{bal9}
\left| 10^{\ell}\alpha^{-n_2}\left(\dfrac{d}{9a(1+\alpha^{n_1-n_2})}\right)-1\right|< \dfrac{3}{\alpha^{n_2-n_3}}.
\end{eqnarray}
We put
\begin{eqnarray*}
\Lambda_2:= 10^{\ell}\alpha^{-n_2}\left(\dfrac{d}{9a(1+\alpha^{n_1-n_2})}\right)-1.
\end{eqnarray*}
As before, in order to apply Theorem \ref{Matveev11} we need to check that $ \Lambda_2\ne 0 $. Suppose that $ \Lambda_2=0 $, then we have
\begin{eqnarray}\label{gal2}
a(\alpha^{n_1}+\alpha^{n_2}) = \dfrac{10^{\ell}\cdot d}{9}.
\end{eqnarray}
To see that this is not true, we again  consider the $ \mathbb{Q} $-automorphism $ \sigma $ of the Galois extension $ \mathbb{Q}(\alpha, \beta) $ over $ \mathbb{Q} $ given by $ \sigma(\alpha):=\beta $ and $ \sigma(\beta):=\alpha $. Now, if $ \Lambda_2=0 $, then $ \sigma(\Lambda_2)=0 $.
Thus, conjugating the relation \eqref{gal2}  under $ \sigma $, and taking absolute values on both sides, we get 
\begin{eqnarray*}
\dfrac{10^{\ell}\cdot d}{9} = |\sigma (a(\alpha^{n_1}+\alpha^{n_2}))|=|b|(|\beta|^{n_1}+|\beta|^{n_2})<  2|b|< \dfrac{2}{3},
\end{eqnarray*}
which is false for $ \ell \ge 2 $ and $ d \ge 1 $. Therefore, $ \Lambda_2\ne 0 $.

So we apply Theorem \ref{Matveev11} with the data
\begin{eqnarray*}
&t:=3, \quad \eta_1:=10, \quad  \eta_2:=\alpha, \quad \eta_3:=\dfrac{d}{9a(1+\alpha^{n_1-n_2})},\quad  b_1:=\ell, \quad b_2:=-n_2, \quad b_3:=1.
\end{eqnarray*}
Since $ \eta_1, \eta_2, \eta_2\in  \mathbb{Q}(\alpha) $, we take the field $ \mathbb{K}:= \mathbb{Q}(\alpha) $ with degree $ D_{\mathbb{K}}:=3 $. Since $ \max\{1, \ell, n_2\} \le n_1$, we take $ B:=n_1 $. Further,
\begin{eqnarray*}
h(\eta_3)&\le& h(d)+h(9) + h(a) + h(1+\alpha^{n_1-n_2})\\&\le& 4\log 3+\dfrac{1}{3}\log 23 +(n_1-n_2)\log \alpha + \log 2\\&<& 1.77\times 10^{14}\log n_1 \quad (\text{by \eqref{bala1}}).
\end{eqnarray*}
Thus, we can take $ A_1:= 3\log 10 $, $ A_2:=\log \alpha $ and $ A_3:=5.31\times 10^{14}\log n_1 $. So, Theorem \ref{Matveev11} tells us that the left-hand side of \eqref{bal9} is bounded below by
\begin{eqnarray*}
\log |\Lambda_2|&>&-1.4 \times 30^{6}\times 3^{4.5}\times 3^{2}(1+\log 3)(1+\log n_1)(3\log 10)(\log \alpha)(5.31\times 10^{14}\log n_1)\\
&>&-1.98\times 10^{28}(\log n_1)^2 \log\alpha.
\end{eqnarray*}
By comparing the above inequality with the right-hand side of \eqref{bal9} we get that 
\begin{eqnarray}\label{bala2}
n_2-n_3 \le 2\times 10^{28}(\log n_1)^2.
\end{eqnarray}
\subsection{Case 3}
 We have that
\begin{eqnarray*}
a\alpha^{n_1}+e(n_1)+a\alpha^{n_2}+e(n_2) + a\alpha^{n_3}+e(n_3) -\dfrac{d\cdot 10^{\ell}}{9} = -\dfrac{d}{9}.
\end{eqnarray*}
This is equivalent to 
\begin{eqnarray*}
a(\alpha^{n_1} +\alpha^{n_2}+\alpha^{n_3})- \dfrac{d\cdot 10^{\ell}}{9} = -\dfrac{d}{9}- e(n_1)-e(n_2)-e(n_3).
\end{eqnarray*}
Thus, we have 
\begin{eqnarray*}
\left|a(\alpha^{n_1} +\alpha^{n_2}+\alpha^{n_3})- \dfrac{d\cdot 10^{\ell}}{9}\right| &\le  &\dfrac{d}{9}+|e(n_1)|+|e(n_2)|+|e(n_3)|\\
&<& 1+3\alpha^{-n_3/2}< 3,
\end{eqnarray*}
and so
\begin{eqnarray}\label{bal10}
\left|a(\alpha^{n_1} +\alpha^{n_2}+\alpha^{n_3})- \dfrac{d\cdot 10^{\ell}}{9}\right|  < 3.
\end{eqnarray}
We divide through \eqref{bal5} by $ a(\alpha^{n_1} +\alpha^{n_2}+\alpha^{n_3})$ to get
\begin{eqnarray*}
\left| 10^{\ell}\alpha^{-n_3}\left(\dfrac{d}{9a(1+\alpha^{n_1-n_3}+\alpha^{n_2-n_3})}\right)-1\right| &<& \dfrac{3\alpha^{-n_1}}{(1+\alpha^{n_2-n_1}+\alpha^{n_3-n_1})}. 
\end{eqnarray*}
Thus, we have
\begin{eqnarray}\label{bal11}
\left| 10^{\ell}\alpha^{-n_3}\left(\dfrac{d}{9a(1+\alpha^{n_1-n_3}+\alpha^{n_2-n_3})}\right)-1\right|< \dfrac{5}{\alpha^{n_1}}.
\end{eqnarray}
We put
\begin{eqnarray*}
\Lambda_3:= 10^{\ell}\alpha^{-n_3}\left(\dfrac{d}{9a(1+\alpha^{n_1-n_3}+\alpha^{n_2-n_3})}\right)-1.
\end{eqnarray*}
As before, in order to apply Theorem \ref{Matveev11} we need to check that $ \Lambda_3\ne 0 $. Suppose that $ \Lambda_3=0 $, then we have
\begin{eqnarray}\label{gal4}
a(\alpha^{n_1}+\alpha^{n_2}+\alpha^{n_3}) = \dfrac{10^{\ell}\cdot d}{9}.
\end{eqnarray}
To see that this is not true, we again  consider the $ \mathbb{Q} $-automorphism $ \sigma $ of the Galois extension $ \mathbb{Q}(\alpha, \beta) $ over $ \mathbb{Q} $ given by $ \sigma(\alpha):=\beta $ and $ \sigma(\beta):=\alpha $. Now, if $ \Lambda_3=0 $, then $ \sigma(\Lambda_3)=0 $.
Thus, conjugating the relation \eqref{gal4}  under $ \sigma $, and taking absolute values on both sides, we get 
\begin{eqnarray*}
\dfrac{10^{\ell}\cdot d}{9} = |\sigma (a(\alpha^{n_1}+\alpha^{n_2}+\alpha^{n_3}))|=|b|(|\beta|^{n_1}+|\beta|^{n_2}+|\beta|^{n_3})<  3|b|< 1,
\end{eqnarray*}
which is false for $ \ell \ge 2 $ and $ d \ge 1 $. Therefore, $ \Lambda_2\ne 0 $.

So we apply Theorem \ref{Matveev11} with the data
\begin{eqnarray*}
&t:=3, \quad \eta_1:=10, \quad  \eta_2:=\alpha, \quad \eta_3:=\dfrac{d}{9a(1+\alpha^{n_1-n_3}+\alpha^{n_2-n_3})},\\ & b_1:=\ell, \quad b_2:=-n_3, \quad b_3:=1.
\end{eqnarray*}
Since $ \eta_1, \eta_2, \eta_2\in  \mathbb{Q}(\alpha) $, we take the field $ \mathbb{K}:= \mathbb{Q}(\alpha) $ with degree $ D_{\mathbb{K}}:=3 $. Since $ \max\{1, \ell, n_3\} \le n_1$, we take $ B:=n_1 $. Further,
\begin{eqnarray*}
h(\eta_3)&\le& h(d)+h(9) + h(a) + h(1+\alpha^{n_1-n_3}+\alpha^{n_2-n_3})\\&\le& 4\log 3+\dfrac{1}{3}\log 23 +((n_1-n_3)+ (n_2-n_3))\log \alpha + 2\log 2\\
&<&6\log 3+ ((n_1-n_2)+2(n_2-n_3))\log\alpha\\
&<& 1.72\times 10^{28}(\log n_1)^2 \quad (\text{by \eqref{bala1} and \eqref{bala2}}).
\end{eqnarray*}
Thus, we can take $ A_1:= 3\log 10 $, $ A_2:=\log \alpha $ and $ A_3:=5.16\times 10^{28}(\log n_1)^2 $. So, Theorem \ref{Matveev11} tells us that the left-hand side of \eqref{bal11} is bounded below by
\begin{eqnarray*}
\log |\Lambda_2|&>&-1.4 \times 30^{6}\times 3^{4.5}\times 3^{2}(1+\log 3)(1+\log n_1)(3\log 10)\\&&\times (\log \alpha)(5.16\times 10^{28}(\log n_1)^2)\\
&>&-1.92\times 10^{42}(\log n_1)^3 \log\alpha.
\end{eqnarray*}
By comparing the above inequality with the right-hand side of \eqref{bal9} we get that 
\begin{eqnarray}\label{bala3}
n_1\le 1.94\times 10^{42}(\log n_1)^3.
\end{eqnarray}
Now, we apply Lemma \ref{gl} on the above inequality  \eqref{bala3} with the data: $ r:=3, ~ H:=1.94\times 10^{42},1 ~ L:=n_1$. We obtain that $ n_1 < 2.7\times 10^{48} $. We record what we have proved
\begin{lemma}\label{bounds}
Let $(N, n_1, n_2, n_3, d, \ell ) $ be the nonnegative integer solutions  to the  Diophantine equation \eqref{Problem}  with $ n_1\ge n_2 \ge n_3 \geq 0 $, $ 1\le d\le 9 $ and $ \ell \ge 2 $. Then we have
\begin{eqnarray*}
\ell < n_1 < 3\times 10^{48}.
\end{eqnarray*}
\end{lemma}

\section{Reducing the bounds}
The bounds ontained in Lemma \ref{bounds} are too large to carry out meaningful computations on the computer. Thus, we need to reduce these bounds. To do so, we return to \eqref{bal6}, \eqref{bal9} and \eqref{bal11} and apply Lemma \ref{Dujjella} via the following procedure.

First, we put
$$
\Gamma_1:=\ell \log 10 - n_1\log \alpha +\log\left(\frac{d}{9a}\right), \quad 1\le d\le 9.
$$
For technical reasons, we assume that $ n_1-n_2 \ge 20 $ and go  to \eqref{bal6}. Note that $ e^{\Gamma_1}-1= \Lambda_1 \neq 0 $. Thus, $ \Gamma_1\neq 0 $. If $ \Gamma_1< 0 $ then
\begin{eqnarray*}
0< |\Gamma_1| < e^{|\Gamma_1|}-1 = |\Lambda_1| < \dfrac{5}{\alpha^{n_1-n_2}}.
\end{eqnarray*}
If $ \Gamma_1> 0 $ then we have that $ |e^{\Gamma_1}-1| <1/2  $. Hence $ e^{\Gamma_1}<2 $. Thus, we get that
\begin{eqnarray*}
0< \Gamma_1 < e^{\Gamma_1}-1 = e^{\Gamma_1}|\Lambda_1| < \dfrac{10}{\alpha^{n_1-n_2}}.
\end{eqnarray*}
Therefore, in both cases, we have that
\begin{eqnarray*}
0< |\Gamma_1|=\left|\ell \log 10 - n_1\log \alpha +\log\left(\frac{d}{9a}\right)\right| < \dfrac{10}{\alpha^{n_1-n_2}}
\end{eqnarray*}
Dividing through the above inequality by $ \log \alpha $, we get
\begin{eqnarray}\label{kapa}
0< \left|\ell\dfrac{\log 10}{\log \alpha} - n_1 + \dfrac{\log (d/(9a))}{\log\alpha}\right| < \dfrac{10}{\alpha^{n_1-n_2}\log\alpha}
\end{eqnarray}
If we put
\begin{eqnarray*}
\tau:=\dfrac{\log 10}{\log \alpha} \quad \text{and} \quad \mu_d:=\dfrac{\log (d/(9a))}{\log\alpha}, \quad 1\le d\le 9,
\end{eqnarray*}
we can rewrite \eqref{kapa} as 
\begin{eqnarray}\label{kapa1}
0< \left|\ell\tau - n_1 + \mu_d\right| < 36\cdot\alpha^{-(n_1-n_2)}
\end{eqnarray}
We now apply Lemma \ref{Dujjella} on \eqref{kapa1}. We put $ M:=3\times 10^{48} $. A quick computer search in \textit{Mathematica} reveals that the convergent
\begin{eqnarray*}
\dfrac{p_{106}}{q_{106}}=\dfrac{177652856036642165557187989663314255133456297895465}{21695574963444524513646677911090250505443859600601}
\end{eqnarray*}
of $ \tau $ is such that $ q_{106}> 6M $ and $ \varepsilon_d \ge 0.0129487>0 $. Therefore, with $ A:=36 $ and $ B:=\gamma $ we calculated each value of $ \log (36q_{106}/\varepsilon_d)/\log\alpha $ and found that all of them are at most $ 432 $. Thus, we have that 
$
n_1-n_2 \le 432.
$

Next, we put
$$
\Gamma_2:=\ell \log 10 - n_2\log \alpha +\log\left(\frac{d}{9a(1+\alpha^{n_1-n_2})}\right), \quad 1\le d\le 9.
$$
For technical reasons, as before we assume that $ n_2-n_3 \ge 20 $ and go  to \eqref{bal9}. Note that $ e^{\Gamma_2}-1= \Lambda_2 \neq 0 $. Thus, $ \Gamma_2\neq 0 $. If $ \Gamma_2< 0 $ then
\begin{eqnarray*}
0< |\Gamma_2| < e^{|\Gamma_2|}-1 = |\Lambda_2| < \dfrac{3}{\alpha^{n_2-n_3}}.
\end{eqnarray*}
If $ \Gamma_2> 0 $ then we have that $ |e^{\Gamma_2}-1| <1/2  $. Hence $ e^{\Gamma_2}<2 $. Thus, we get that
\begin{eqnarray*}
0< \Gamma_2 < e^{\Gamma_2}-1 = e^{\Gamma_2}|\Lambda_2| < \dfrac{6}{\alpha^{n_2-n_3}}.
\end{eqnarray*}
Therefore, in both cases, we have that
\begin{eqnarray*}
0< |\Gamma_2|=\left|\ell \log 10 - n_1\log \alpha +\log\left(\frac{d}{9a(1+\alpha^{n_1-n_2})}\right)\right| < \dfrac{6}{\alpha^{n_2-n_3}}
\end{eqnarray*}
Dividing through the above inequality by $ \log \alpha $, we get
\begin{eqnarray}\label{kapaa}
0< \left|\ell\dfrac{\log 10}{\log \alpha} - n_2 + \dfrac{\log (d/(9a(1+\alpha^{n_1-n_2})))}{\log\alpha}\right| < \dfrac{6}{\alpha^{n_2-n_3}\log\alpha}
\end{eqnarray}
We put
\begin{eqnarray*}
\tau:=\dfrac{\log 10}{\log \alpha} \quad \text{and} \quad \mu_{d,k}:=\dfrac{\log (d/(9a(1+\alpha^{k}))}{\log\alpha}, \quad 1\le d\le 9, \quad 1\le k \le 432,
\end{eqnarray*}
 where $ k:=n_1-n_2 $. We can rewrite \eqref{kapaa} as 
\begin{eqnarray}\label{kapaa1}
0< \left|\ell\tau - n_2 + \mu_{d,k}\right| < 22\cdot\alpha^{-(n_2-n_3)}
\end{eqnarray}
We now apply Lemma \ref{Dujjella} on \eqref{kapaa1}. We put $ M:=3\times 10^{48} $. A quick computer search in \textit{Mathematica} reveals that the $ 106 $-th convergent
of $ \tau $ is such that $ q_{106}> 6M $ and $ \varepsilon_{d,k} \ge 0.000134829>0 $ for all $ 1\le d\le 9 $ and $ 1\le k \le 432 $ except for the case $ \varepsilon_{9,11} $, which is always negative. Therefore, with $ A:=22 $ and $ B:=\alpha $ we calculated each value of $ \log (22q_{106}/\varepsilon_{d,k})/\log\alpha $ and found that all of them are at most $ 446 $. Thus, we have that 
$
n_2-n_3 \le 446.
$

The problem in the case of $ \varepsilon_{9,11} $ is due to the fact that
\begin{eqnarray}\label{kap}
\dfrac{1}{\alpha^{9}}= \dfrac{3\alpha^{2}-1}{\alpha(\alpha+1)(\alpha^{11}+1)}
\end{eqnarray}
Thus, if we consider the identity \eqref{kap}, the inequality \eqref{kapaa} becomes
\begin{eqnarray}\label{kapaax}
0< \left|\tau - \dfrac{(n_2+9}{\ell}\right| < \dfrac{6}{\alpha^{n_2-n_3}\ell\log\alpha}.
\end{eqnarray}
In this case we apply the classical result from Diophantine approximation given in Lemma \ref{Legendre}. We assume that $ n_2-n_3 $ is so large that the right-hand side of the inequality in  \eqref{kapaax} is smaller than $ 1/(2\ell^2) $. This certainly holds if
\begin{eqnarray}\label{DNA1}
\alpha^{n_2-n_3}>12\ell/\log\alpha.
\end{eqnarray}
Since $ \ell< n_1< 3\times 10^{48} $, it follows that the last inequality \eqref{DNA1} holds provided that $ n_2-n_3\ge 415 $, which we now assume. In this case $ r/s:=(n_2+9)/\ell $ is a convergent of the continued fraction of $ \tau: = \log 10 / \log \alpha $ and $ \ell < 3\times 10^{48} $. We are now set to apply Lemma \ref{Legendre}. 

We write $ \tau: =[a_0; a_1, a_2, a_3, \ldots] = [8; 5, 3, 3, 1, 5, 1, 8, 4, 6, 1, 4, 1, 1, 1, 9, 1, 4, 4, 9, \ldots] $ for the continued fraction of $ \tau $ and $ p_k/q_k $ for the $ k- $th convergent. We get that $ r/s=p_j/q_j $ for some $ j\le 106  $. Furthermore, putting $ a(M):=\max\{a_j: j=0,1, \ldots, 106\} $, we get $ a(M):=564 $. By Lemma \ref{Legendre}, we get
\begin{eqnarray*}
\dfrac{1}{566\ell^2}=\dfrac{1}{(a(M)+2)\ell^{2}}\le \left|\tau-\dfrac{r}{s}\right|<  \dfrac{6}{\alpha^{n_2-n_3}\ell\log\alpha},
\end{eqnarray*}
giving
\begin{eqnarray*}
\alpha^{n_2-n_3}< \dfrac{566\times 6\ell }{\log\alpha} < \dfrac{566\times 6 \times 3\times 10^{48}}{\log\alpha},
\end{eqnarray*}
leading to $ n_2-n_3\le 435 $. Thus, in both cases we have that $ n_2-n_3\le 446 $.

Lastly, we put 
$$
\Gamma_3:=\ell \log 10 - n_3\log \alpha +\log\left(\frac{d}{9a(1+\alpha^{n_1-n_3}+\alpha^{n_2-n_3})}\right), \quad 1\le d\le 9.
$$
For technical reasons, as before we assume that $ n_1 \ge 20 $ and go  to \eqref{bal11}. Note that $ e^{\Gamma_3}-1= \Lambda_3 \neq 0 $. Thus, $ \Gamma_3\neq 0 $. If $ \Gamma_3< 0 $ then
\begin{eqnarray*}
0< |\Gamma_3| < e^{|\Gamma_3|}-1 = |\Lambda_3| < \dfrac{5}{\alpha^{n_1}}.
\end{eqnarray*}
If $ \Gamma_3> 0 $ then we have that $ |e^{\Gamma_3}-1| <1/2  $. Hence $ e^{\Gamma_3}<2 $. Thus, we get that
\begin{eqnarray*}
0< \Gamma_3 < e^{\Gamma_3}-1 = e^{\Gamma_3}|\Lambda_3| < \dfrac{10}{\alpha^{n_1}}.
\end{eqnarray*}
Therefore, in both cases, we have that
\begin{eqnarray*}
0< |\Gamma_3|=\left|\ell \log 10 - n_3\log \alpha +\log\left(\frac{d}{9a(1+\alpha^{n_1-n_2}+\alpha^{n_2-n_3})}\right)\right| < \dfrac{10}{\alpha^{n_1}}
\end{eqnarray*}
Dividing through the above inequality by $ \log \alpha $, we get
\begin{eqnarray}\label{kapaaa}
0< \left|\ell\dfrac{\log 10}{\log \alpha} - n_3 + \dfrac{\log (d/(9a(1+\alpha^{n_1-n_3}+\alpha^{n_2-n_3}))}{\log\alpha}\right| < \dfrac{10}{\alpha^{n_1}\log\alpha}
\end{eqnarray}
We put
\begin{eqnarray*}
\tau:=\dfrac{\log 10}{\log \alpha} \quad \text{and} \quad \mu_{d,k,s}:=\dfrac{\log \left(d/(9a(1+\alpha^{k}+\alpha^{s}))\right)}{\log\alpha}, \quad 1\le d\le 9,
\end{eqnarray*}
 where $ 1\le k:=n_1-n_3 =(n_1-n_2)+(n_2-n_3) \le 878 $ and $ 1\le s:=n_2-n_3 \le 446 $. We can rewrite \eqref{kapaaa} as 
\begin{eqnarray}\label{kapaaa1}
0< \left|\ell\tau - n_3 + \mu_{d,k,s}\right| < 36\cdot\alpha^{-n_1}
\end{eqnarray}
We now apply Lemma \ref{Dujjella} on \eqref{kapaaa1}. We put $ M:=3\times 10^{48} $. A quick computer search in \textit{Mathematica} reveals that the $ 106 $-th convergent
of $ \tau $ is such that $ q_{106}> 6M $ and $ \varepsilon_{d,k,s} \ge 0.000125>0 $. Therefore, with $ A:=36 $ and $ B:=\alpha $ we calculated each value of $ \log (36q_{106}/\varepsilon_{d,k,s})/\log\alpha $ and found that all of them are at most $ 485 $. Thus, we have that 
$
n_1 \le 485.
$
This contradicts our assumption that $ n_1> 500 $. Hence,  Theorem \ref{Main} holds. \qed

\begin{acknowledgements}
The author was supported by the Austrian Science Fund (FWF) grants: F5510-N26 -- Part of the special research program (SFB), ``Quasi Monte Carlo Metods: Theory and Applications'' and W1230 --``Doctoral Program Discrete Mathematics''.
\end{acknowledgements}

\def\cprime{$'$}

\end{document}